\magnification=1200
\overfullrule=0pt
\centerline {\bf A conjecture implying the existence of
non-convex Chebyshev sets}\par
\centerline {\bf in infinite-dimensional Hilbert spaces}\par
\bigskip
\bigskip
\centerline {BIAGIO RICCERI}\par
\bigskip
\bigskip
{\it Abstract.} In this paper, we propose the study of a conjecture whose
positive solution would provide an example of a non-convex Chebyshev set
in an infinite-dimensional real Hilbert space.\par
\bigskip
{\it Keywords:} Chebyshev set, Hilbert space, convexity.\par
\bigskip
{\it AMS 2000 Subject Classification:} 41A50, 41A65.
\bigskip
\bigskip
\bigskip
\bigskip
Here and in the sequel, $(X,\langle\cdot,\cdot\rangle)$ is a separable real Hilbert space, with norm
$\|\cdot\|$. A non-empty set $C\subset X$ is said to be a 
Chebyshev set if, for each $x\in X$, there exists a unique $y\in C$ such that
$$\|x-y\|=\inf_{z\in C}\|x-z\|\ .$$
Clearly, each closed convex set is a Chebyshev one. A natural question is:
must any Chebyshev $C\subset X$ be convex ? We refer to
the surveys [1], [5] for
a thorough discussion of the subject. In particular,
it is well-known that any sequentially weakly closed Chebyshev set
$C\subset X$ is convex. Hence, if $X$ is finite-dimensional, the
answer to the above question is "yes".
\smallskip
However, since [7], it is a quite common feeling that if $X$ is
infinite-dimensional,
then $X$ contains some non-convex Chebyshev set (see also [6] for a recent
contribution in this direction). Maybe, this is the most
important conjecture in best approximation theory. \par
\smallskip
A much more recent (and less known) problem is: if $f:X\to {\bf R}$ is
a lower semicontinuous function
 such that, for
each $y\in X$ and each $\lambda>0$, the function $x\to \|x-y\|^2+\lambda f(x)$
has a unique global minimum, must $f$ be convex ? For this problem too, 
 the answer is "yes" if $X$ is finite-dimensional ([11], Corollary 3.8). See also  
Corollary 5.2 of [2]
for another partial answer.\par
\smallskip
The aim of the present paper is to show that if the second problem
has a qualified negative answer, then the same happens for the first one.\par
\smallskip
In the sequel, 
$L^2([0,1],X)$ is the usual space of all (equivalence classes of)
 measurable functions $u:[0,1]\to
X$ such that $\int_0^1 \|u(t)\|^2dt<+\infty$, endowed with
the scalar product
$$\langle u,v\rangle_{L^2_X}=\int_0^1\langle u(t),v(t)\rangle dt$$
The norm induced by $\langle\cdot,\cdot\rangle_{L^2_X}$ is denoted
by $\|\cdot\|_{L^2_X}$.\par
\smallskip
Let us start with the following\par
\medskip
DEFINITION 1. - Let $Y$ be a non-empty set and ${\cal F}$ a family of
subsets of $Y$.\par
We say that ${\cal F}$ has the compactness-like property if
every subfamily of ${\cal F}$ satisfying the finite intersection property has
a non-empty intersection. \par
\medskip
We have the following characterization which is due to C. Costantini ([3]):\par
\medskip
PROPOSITION 1. - {\it Let $Y$ be a non-empty set, let
${\cal F}$ be a family of subsets of $Y$ and let $\tau$ be the topology
on $Y$ generated by the family $\{Y\setminus C\}_{C\in {\cal F}}$.\par
Then, the following assertions are equivalent:\par
\noindent
$(i)$\hskip 5pt Each member of ${\cal F}$ is $\tau$-compact.\par
\noindent
$(ii)$\hskip 5pt The family ${\cal F}$ has the compactness-like property.\par
\noindent
$(iii)$\hskip 5pt The space $Y$ is $\tau$-compact.}\par
\medskip
We then formulate the following \par
\medskip
CONJECTURE 1. - If $X$ is infinite-dimensional, there
exist a non-convex Borel function $f:X\to {\bf R}$,
$r\in ]\inf_Xf,\sup_X f[$ and $\gamma\in ]0,+\infty]$, with
the following properties:\par
\noindent
$(a)$
$$\sup_{x\in X}{{|f(x)|}\over {1+\|x\|^2}}<+\infty\ ;$$
\noindent
$(b)$\hskip 5pt 
for each $y\in X$ and each
$\lambda\in ]0,\gamma[$, the function
$$x\to \|x-y\|^2+\lambda f(x)$$
has a unique global minimum in $X$, say
$\hat x_{y,\lambda}$; moreover, the map
$y\to \hat x_{y,\lambda}$ is Borel and one has
$$\|\hat x_{y,\lambda}\|\leq c_{\lambda}(1+\|y\|)$$
where $c_\lambda$ is independent of $y$\ ;\par
\noindent
(c) \hskip 5pt if $\gamma<+\infty$,
for each $y\in f^{-1}(]r,+\infty[)$, the function
$$x\to \|x-y\|^2+ \gamma f(x)$$
has no global minima in $X$\ ;\par
\noindent
$(d)$\hskip 5pt
 for each $v\in L^2([0,1],X)$, with
$\int_0^1f(v(t))dt>r$, the family
$$\left \{\left \{ u\in L^2([0,1],X) :
\int_0^1\|u(t)-v(t)\|^2dt+\lambda \int_0^1 f(u(t))dt\leq \rho\right \} :
\lambda\in ]0,\gamma[, \rho\in {\bf R}\right \} $$
has the compactness-like property.\par
\medskip
Our result reads as follows:\par
\medskip
THEOREM 1. - {\it Assume that Conjecture 1 is true and let
$f$ be a function satisfying it.\par
Then,
$$\left \{ u\in L^2([0,1],X): \int_0^1 f(u(t))dt\leq r\right \}$$
 is a non-convex Chebyshev set.}\par
\medskip
To prove Theorem 1, we need the following two results. \par
\medskip
THEOREM A. - {\it Let $Y$ be a non-empty set, $\eta\in ]0,+\infty]$ and $\varphi, \psi:
Y\to {\bf R}$ two functions such that
 the function $\varphi+\lambda\psi$
has a unique global minimum if $\lambda\in [0,\eta[$, while has no
global minima if $\eta<+\infty$ and $\lambda=\eta$.
Moreover, if $y_0$ is the only global
minimum of $\varphi$, assume that $\inf_Y\psi<\psi(y_0)$. Finally,
assume that the family 
$$\{\{y\in Y : \varphi(y)+\lambda\psi(y)\leq \rho\} : \lambda\in ]0,\eta[, \rho\in {\bf R}\}$$
has the compactness-like property.\par
Then, for each $\rho\in ]\inf_Y\psi,\psi(y_0)[$, the restriction
of the function $\varphi$ to $\psi^{-1}(\rho)$ has a unique global minimum.}\par
\medskip
THEOREM B. - {\it Let $f:X\to {\bf R}$ be a Borel function such
that
$$\sup_{x\in X}{{|f(x)|}\over {1+\|x\|^2}}<+\infty\ .$$
Assume that, for some $\rho\in ]\inf_Xf,\sup_Xf[$, the
set 
$$\left \{ u\in L^2([0,1],X): \int_0^1 f(u(t))dt\leq \rho\right \}$$
 is weakly closed.\par
Then, $f$ is convex.}\par
\medskip
Theorem A, via Proposition 1, is a 
direct consequence of a variant of Theorem 1 of [9] (see also the proof
of Theorem 1 of [10]), 
while
Theorem B has been proved by R. Landes in [8]. 
\medskip
{\it Proof of Theorem 1}. 
Fix $\lambda\in ]0,\gamma[$, 
  $v\in L^2([0,1],X)$, with $\int_0^1 f(v(t))dt>r$, and put
$$\omega_{v,\lambda}(t)=\hat x_{v(t),\lambda}$$
for all $t\in [0,1]$.  From $(a)$ and $(b)$,
it clearly follows that the function $\omega_{v,\lambda}$
 belongs to $L^2([0,1],X)$. 
If $u\in L^2([0,1],X)$ and $u\neq \omega_{v,\lambda}$, we have
$$\|\omega_{v,\lambda}(t)-v(t)\|^2+\lambda f(\omega_{v,\lambda}(t))\leq  
\|u(t)-v(t)\|^2+\lambda f(u(t))$$
for all $t\in [0,1]$, the inequality being strict in a subset of $[0,1]$ with
positive measure. Then, by integrating, we get
$$\int_0^1\|\omega_{v,\lambda}(t)-v(t)\|^2dt+\int_0^1\lambda f(\omega_{v,\lambda}(t))dt< 
\int_0^1\|u(t)-v(t)\|^2dt+\lambda \int_0^1f(u(t))dt\ .$$
Therefore, $\omega_{v,\lambda}$ is the only global minimum in $L^2([0,1],X)$ of
the functional $$u\to \int_0^1 \|u(t)-v(t)\|^2dt+\lambda\int_0^1f(u(t))dt\ .$$
 Now, assume that $\gamma<+\infty$. Put
$$A_v=\{t\in [0,1] : f(v(t))>r\}\ .$$
Since $\int_0^1f(v(t)dt>r$, the measure of $A_v$ is positive.
 We show that the functional
$$u\to \int_0^1 \|u(t)-v(t)\|^2dt+\gamma \int_0^1f(u(t))dt$$
has no global minima in $L^2([0,1],X)$. Indeed, fix
$u\in L^2([0,1],X)$. It is easy to check that the function
$(t,x)\to \|x-v(t)\|^2+\gamma f(x)$ is
 ${\cal L}([0,1])\otimes {\cal B}(X)$-measurable,
where ${\cal L}([0,1])$ and ${\cal B}(X)$ denote the
Lebesgue and the Borel $\sigma$-algebras of subsets of $[0,1]$ and
$X$, respectively. So, by Theorem 2.6.40 of [4], the function
 $t\to \inf_{x\in X}(\|x-v(t)\|^2+f(x))$
is measurable. On the other hand, in view of $(c)$, we have
$$\inf_{x\in X}(\|x-v(t)\|^2+\gamma f(x))<\|u(t)-v(t)\|^2+\gamma f(u(t))$$
for all $t\in A_v$. Consequently, 
 we can apply Theorem 4.3.7 of [4] to get 
a measurable function $\xi:[0,1]\to X$ such that
$$\|\xi(t)-v(t)\|^2+\gamma f(\xi(t))<\|u(t)-v(t)\|^2+\gamma f(u(t))$$
for all $t\in A_v$. Finally, choose a set $B\subset A$ with positive
measure such that $\xi$ is bounded in $B$ and put
$$w(t)=\cases {\xi(t) & if $t\in B$\cr & \cr u(t) & if $t\in [0,1]\setminus B\ .$\cr}$$
Clearly, $w\in L^2([0,1],X)$ and one has
$$\int_0^1 \|w(t)-v(t)\|^2dt+\gamma
\int_0^1 f(w(t))dt<\int_0^1 \|u(t)-v(t)\|^2dt+\gamma\int_0^1 f(u(t))dt$$
which proves our claim. 
At this point, we can apply Theorem A taking
$$Y=L^2([0,1],X)\ ,$$
$$\eta=\gamma\ ,$$
$$\varphi(u)=\|u-v\|_{L^2_X}^2$$
and
$$\psi(u)=\int_0^1f(u(t))dt\ .$$
Then, there exists a unique $u\in \psi^{-1}(r)$ such that
$$\|v-u\|_{L^2_X}=\hbox {\rm dist}(v,\psi^{-1}(r))\ .$$
We now claim that such an $u$ is the unique
point of $\psi^{-1}(]-\infty,r])$ such that
$$\|v-u\|_{L^2_X}=\hbox {\rm dist}(v,\psi^{-1}(]-\infty,r]))\ .$$
This amounts to show that if $w\in \psi^{-1}(]-\infty,r])$ is such that
$$\|v-w\|_{L^2_X}=\hbox {\rm dist}(v,\psi^{-1}(]-\infty,r]))\ ,\eqno{(1)}$$
then $\psi(w)=r$. Arguing by contradiction, assume that
$\psi(w)<r$. For each measurable set $A\subset [0,1]$, put
$$h_A(t)=\cases {v(t) & if $t\in A$\cr & \cr w(t) & if $t\in [0,1]\setminus A\ .$\cr}$$
Also, set
$$D=\{h_A : A\subset [0,1],\hskip 3pt A\hskip 3pt \hbox {\rm measurable}\}\ .$$
It is not hard to check that $D$ is decomposable ([4],
p. 452). Moreover, it is clear that
$v, w\in D$ and that 
$$\|v-h\|<\|v-w\| \eqno{(2)}$$ for all $h\in D\setminus \{v,w\}$. By
Corollary
4.5.13 of [4],
the set $\psi(D)$ is an interval. Consequently, there exists $h\in D\setminus \{v,w\}$
such that $\psi(h)=r$. This implies a contradiction, in view of $(1)$ and $(2)$.
So, $\psi^{-1}(]-\infty,r])$ is a Chebyshev set in
$L^2([0,1],X)$. Finally, this set is not convex. Indeed, if it was convex, being closed,
it would be weakly closed. Then, by Theorem B, the function
$f$ would be convex, against the assumptions.
\hfill $\bigtriangleup$

\vfill\eject
\centerline {\bf References}\par
\bigskip
\bigskip
\noindent
[1]\hskip 5pt V. S. BALAGANSKII and L. P. VLASOV, {\it The problem
of the convexity of Chebyshev sets}, Russian Math. Surveys, {\bf 51}
(1996), 1127-1190.\par
\smallskip
\noindent
[2]\hskip 5pt F. BERNARD and L. THIBAULT, {\it Prox-regular
functions in Hilbert spaces}, J. Math. Anal. Appl., {\bf 303} (2005),
1-14.\par
\smallskip
\noindent
[3]\hskip 5pt C. COSTANTINI, Personal communication.\par
\smallskip
\noindent
[4]\hskip 5pt Z. DENKOWSKI, S. MIG\'ORSKI and N. S. PAPAGEORGIOU,
{\it An Introduction to Nonlinear Analysis: Theory}, Kluwer Academic
Publishers, 2003.\par
\smallskip
\noindent
[5]\hskip 5pt F. DEUTSCH, {\it The convexity of Chebyshev sets in Hilbert
space}, in {\it Topics in polynomials of one and several variables and their
applications}, 143-150, World Sc. Publishing, 1993.\par
\smallskip
\noindent
[6]\hskip 5pt F. FARACI and A. IANNIZZOTTO, {\it Well posed optimization
problems and nonconvex Chebyshev sets in Hilbert spaces}, SIAM J. Optim.,
{\bf 19} (2008), 211-216.\par
\smallskip
\noindent
[7]\hskip 5pt V. L. KLEE, {\it Remarks on nearest points in normed linear spaces},
in Proc. Colloquium on Convexity (Copenhagen, 1965), 168-176, Kobenhavns Univ.
Mat. Inst., Copenhagen, 1967.\par
\smallskip
\noindent
[8]\hskip 5pt R. LANDES, {\it On a necessary condition in the calculus of
variations}, Rend. Circ. Matematico Palermo, {\bf 41} (1992), 369-387.\par
\smallskip
\noindent
[9]\hskip 5pt B. RICCERI, {\it Well-posedness of constrained minimization
problems via saddle-points}, J. Global Optim., {\bf 40} (2008), 389-397.\par
\smallskip
\noindent
[10]\hskip 5pt B. RICCERI, {\it Multiplicity of global minima for
parametrized functions}, Rend. Lincei Mat. Appl., {\bf 21} (2010),
47-57.\par
\smallskip
\noindent
[11]\hskip 5pt X. WANG, {\it On Chebyshev functions and Klee functions},
J. Math. Anal. Appl., {\bf 368} (2010), 293-310.

\bigskip
\bigskip
\bigskip
\bigskip
Department of Mathematics\par
University of Catania\par
Viale A. Doria 6\par
95125 Catania\par
Italy\par
{\it e-mail address}: ricceri@dmi.unict.it

\bye